# Arithmetic on Balanced Parentheses: The case of Ordered Motzkin Words

Gennady Eremin

ergenns@gmail.com

November 5, 2019

**Abstract.** In this article, we establish a total lexicographical order on the set of Motzkin words. Elements are ordered similarly to Natural Numbers in accordance with known rules (axioms). As a result, we were able to obtain arithmetic and logical operations on the elements of the ordered sequence, *Motzkin Row*. This sequence consists of balanced brackets without leading zeros, with the exception of the initial word "0". It is the word "0" as well as the alphabetical symbol "0" that are analogues of numeric zero in the corresponding operations. Logical operations allow you to navigate Motzkin Row. Operations on words are accompanied by index equations, *index polynomials*.

*Key Words*: Motzkin word, lexicographical order, lexicographical sequence, arithmetic on parentheses, index polynomial.

## 1 Introduction

There is a growing interest in ordered non-numerical sets. We note two papers [BP14, Fan19] in which the authors investigate a partial order on well-formed parentheses strings (Dick words and Motzkin words). For example, Pallo and Baril [BP14] determine the distance between the specific Motzkin words in the Tamari lattice.

The author establishes a total order on balanced brackets in accordance with the alphabetical order (see order on Dyck words [Ere19]). In this paper, we work with unique Motzkin words (without leading zeros). The lexicographic sequence is constructed in accordance with the rules (axioms) of the set of natural numbers. As a result, partial operations of addition and subtraction of Motzkin words began to work. Logical operations also appeared.

1.1. **Motzkin words.** In discrete mathematics, *Motzkin words* (named after Theodore Motzkin) are of particular interest. We will treat a Motzkin word as balanced parentheses supplemented by zeros. Let's convert the simple arithmetic expression

$$y \times [(82 - z) / (m + t) - 7] + 359 - x.$$

First, replace the operations and operands with zeros. Next, we change all the brackets to parentheses, keeping the orientation. The result is Motzkin word

(1.1) $\qquad\qquad\qquad$ 00((000)0(000)00)0000.

Two conditions are true for such a string:

- the number of left and right parentheses is the same (a Motzkin word can consist of only zeros, no parentheses, or of only parentheses, no zeros);
- every substring must contain at least as many open parentheses as closed ones.

For any Motzkin word, we can one way to highlight *matched pairs* of bidirectional parentheses. In (1.1) one matched pair is shown in red. Matched pairs must be correctly nested. For the matched pair two conditions are met:
- the left parentheses precedes the right one;
- inside the matched pair, it's either a Motzkin word or nothing (the case of *adjacent parentheses*).

The set of consecutive Motzkin words is also a Motzkin word. Any Motzkin word is created using "building blocks" of the following three types: a left parenthesis, a right parenthesis, and zero. The simplest word contains only zeros. Further, analyzing Motzkin words, we will focus on *leading* (initial) and *final* zeros. For example, the Motzkin word (1.1) has two leading and four final zeros.

In enumerative problems of combinatorics, we usually count the number of elements of some set. There are only two Motzkin words of size 2, *2-word*, namely:

$$00, \ ().$$

The three-character Motzkin word, *3-word*, can be obtained in four variants:

$$000, \ 0(\,), \ (0), \ (\,)0.$$

The first two 3-words (highlighted in green) are inherited from the 2-words by adding one leading zero (concatenation of zero and a 2-word). The last two 3-words are unique. In the set of Motzkin 4-words there are the following nine elements:

$$0000, \ 00(\,), \ 0(0), \ 0(\,)0, \ (00), \ (0)0, \ ((\,)), \ (\,)00, \ (\,)(\,).$$

Similarly, the first four 4-words (again green) are inherited from the 3-words (extra zero at the beginning). The last five 4-words are unique. Naturally, this procedure can be continued. Note that if we denote by $U_n$ the number of the unique *n*-words, then the number of all Motzkin 4-words, $M_4$, is

(1.2) $\quad M_4 = M_3 + U_4 = M_2 + U_3 + U_4 = M_1 + U_2 + U_3 + U_4 = 1+1+2+5 = 9.$

Thus, the set of Motzkin *n*-words can be divided into two disjoint subsets. The first subset combines the inherited elements; these are Motzkin (*n*–1)-words with an extra zero at the beginning of the code. The elements of the second subset are unique, their codes begin with the left parenthesis, and these elements are usually the majority. It is easy to check, if to analyze Motzkin numbers [Slo19, sequence A001006]. Here is the first Motzkin numbers $M_n$ for $n = 0, 1, 2, \ldots$:

(1.3) $\quad$ [1,] 1, 2, 4, 9, 21, 51, 127, 323, 835, 2188, 5798, 15511, 41835 …

We have enclosed the initial integer in brackets because some authors omit the 0th element (see [Wei19] or [Fan19, page 3]). Easy to see, $M_n > 2M_{n-1}$, $n \geq 4$. Usually there are more unique Motzkin *n*-words than inherited ones. But for the first Motzkin numbers we have $M_3 = 2M_2 = 4$ and $M_2 = 2M_1 = 2$. The 2-word "00" is inherited from the Motzkin 1-word "0" and the other 2-word "( )" is unique.

Let's go down another step, here we have $M_1 = M_0 = 1$. What about the 1-word "0"? Is this word inherited or unique? If Motzkin numbers are indexed from 0, then we allow the empty word $\epsilon$ (virtual 0-word). It is logical to consider "0" as the inherited 1-word, since "$0\epsilon$" = "0". But in this paper, we work with the real Motzkin words, we do not use the empty word and number $M_0$. In this regard, the 1-word "0" is considered unique, and this is the only unique Motzkin word with a leading zero.

1.2. **Axioms of the lexicographic sequence.** In [Ere19], the author proposes to establish a lexicographical order using the external features of the sequence of natural numbers. Let's briefly repeat the relevant axioms; these axioms we put in the basis of *Motzkin Row*. The first axiom is obvious: there are no repeats among natural numbers. We can add leading zeros to any integer, but in this case the value of the number will not change. So the first axiom is about the uniqueness of the elements.

**Axiom 1.1** (unique elements). *In the lexicographical sequence, all elements are unique.*

Next, natural numbers are listed in increasing order and, at first glance, the order is determined by the code length. Integers are distributed over *ranges*. Single-digit numbers, 1-*range*, are listed first, followed by double-digit integers, 2-*range*, and so on. Let's call it a *primary order*.

Again, adding leading zeros only formally changes the code length for an integer (the number remains the same). Copies with leading zeros do not modify the set of natural numbers. Leading zeros can be added temporarily for some procedures, for example, summing columns.

The following axiom fixes the primary (external) order.

**Axiom 1.2** (primary order). *In the lexicographic sequence, elements are arranged in increasing order of the code length.*

In ranges, integers are sorted according to the alphabetical order $0 < 1 < 2 < ... < 9$. In the *k*-range, $k > 1$, a minimum element is 10...0 (*k*–1 of zeros), a maximum element is 99…9 (*k* of nines), and the number of elements is $10^k - 10^{k-1} = 9 \times 10^{k-1}$. It is logical to call sorting in ranges an *internal order*.

An alphabet of the Motzkin word is ternary: {0, (, )}; we have to choose a total alphabetical order to sort the elements within the ranges of Motzkin Row. The third axiom establishes the order in ranges.

**Axiom 1.3** (internal order). *In each range of the lexicographic sequence, objects are sorted according to the given total alphabetical order.*

In a natural number, alphabet symbols are *free*, their location is arbitrary. Digits are not related to each other (unlike the brackets in a Motzkin word where each left parenthesis must have a paired right one and vice versa). However, there is some limitation for zero, the smallest symbol of the alphabet. The next Axiom 1.4 is related to the special status of zero in natural numbers.

A natural number does not start with zero. The singular 0 is an exception to the rule (if zero is a natural number). In Motzkin word, zero is also present, and this symbol is located freely within the code. Motzkin words often begin with zero, and something needs to be done. In the meantime, we formulate the last axiom based on the rule of using zero in the code of integers.

**Axiom 1.4** (about minimal free symbol). *In the lexicographic sequence, an element with code length $\geq 2$ cannot begin with a free character that has the least weight in the alphabet.*

## 2 Motzkin Row

We have removed the empty word $\epsilon$ from our analysis (and respectively we got rid of the element $M_0$ in A001006); the axioms described above are generally applicable to construct a lexicographical sequence of Motzkin words. However, we still need a total alphabetical order; this order is natural:

(2.1) $$0 < ( < )$$

In a Motzkin word, any matched pair of brackets starts with the left parenthesis, so it's logical to take the weight of the left parenthesis less than the right one. The symbol "0" resembles zero in integers (it has a minimum weight and is free). It remains for us to decide what to do with leading zeros (except for the 1-word "0").

The order (2.1) allows us to place Motzkin words within a range. In two different words of the same size (without leading zeros), let's look at the characters in pairs from left to right; then the first mismatch will determine the smaller (larger) element in according to (2.1). To determine the next / previous element for a given Motzkin word, we look through the characters from right to left and select the first character that can be replaced with a larger / smaller one. In this case, the tail of the word is selected accordingly.

2.1. **Ranges.** Let's write out the set of Unique Motzkin Words, *Motzkin Row*, in order of increasing size and in accordance with the order (2.1) in the ranges:

$\mathfrak{M} = \{$0, (), (0), ()0, (00), (0)0, (()), ()00, ()(), (000), (00)0, (0()), (0)00, (0)(), ((0)), (()0), (())0, ()000, ()0(), ()(0), ()()0, (0000), (000)0, …$\}$

Elements of the infinite set $\mathfrak{M}$ are indexed starting from zero: $w_0 \equiv 0$, $w_1 \equiv ()$, $w_2 \equiv (0)$, and so on. We identify $w_i$ ($i$ is a specific nonnegative integer) and the $i$-th word of $\mathfrak{M}$. Let's call words of $\mathfrak{M}$ *constants*. For example, constant $(000)0 \equiv w_{22}$.

Let's denote by $U_n$ the cardinality of unique Motzkin words of size $n$, i.e., the length of the *n-range* $\mathfrak{M}_n$ (or $\#\mathfrak{M}_n$). We highlighted in red the first elements in ranges. Obviously (see (1.2)),

$$\mathfrak{M} = \cup_{n \geq 1} \mathfrak{M}_n; \quad \#\mathfrak{M}_1 = U_1 = 1, \quad \#\mathfrak{M}_n = U_n = M_n - M_{n-1}, n \geq 2.$$

The numbers $U_n$ for $n = 1, 2, \ldots$ form the sequence:

(2.2)    1, 1, 2, 5, 12, 30, 76, 196, 512, 1353, 3610, 9713, 26324, 71799, 196938, ….

Each *n*-range, $n \geq 2$, has the following minimum and maximum *n*-words:

$$\min \mathfrak{M}_n = (00\ldots0), n-2 \text{ of zeros, and } \max \mathfrak{M}_n = ()()\ldots()[0], n \geq 2.$$

The maximum element of size *n* ends with zero if *n* is odd, and the number of pairs of adjacent parentheses is $\lfloor n/2 \rfloor$ (we discard half). The maximum and minimum are coincide in singleton ranges $\mathfrak{M}_1 = \{0\}$ and $\mathfrak{M}_2 = \{()\}$.

We construct $\mathfrak{M}$ in the image and likeness of natural numbers

$$\mathbb{N}_0 = \{0, 1, 2, \ldots, 9, 10, \ldots, 99, 100, \ldots, 999, 1000, \ldots, 9999, 10000, \ldots\}.$$

We add that both sets $\mathfrak{M}$ and $\mathbb{N}_0$ have the same 0th elements $w_0 = n_0 = 0$, which have similar properties (and we'll talk about this later). There are other coincidences as well, and we will discuss them in connection with arithmetic and logical operations.

The indices of the elements of $\mathfrak{M}$ are natural numbers. We can say that the set $\mathbb{N}_0$ indexes the set $\mathfrak{M}$, or $\mathbb{N}_0$ is an index set of $\mathfrak{M}$. (The set of natural numbers is self-indexing; the index of the element is equal to its value, i.e., $n_i = i$.) In the ranges of Motzkin Row, minimal and maximum elements are indexed by Motzkin numbers. For any $x \in \mathfrak{M}$, let ind $x$ denote the index of $x$ in $\mathfrak{M}$. Easy to see

(2.3)    ind min $\mathfrak{M}_n = M_{n-1}$   and   ind max $\mathfrak{M}_n = M_n - 1$, $n \geq 2$.

**Definition 2.1**. Index expressions in which Motzkin numbers appear are called *index polynomials*.

In the general case, for $x, y \in \mathfrak{M}$, we will write $x > y$ if ind $x >$ ind $y$. About a thousand elements of Motzkin Row are given in the Addendum, which we will use from time to time.

2.2. **Blocks.**  Recall that in a Motzkin word, all parentheses are divided into matched pairs. A matched pair of parentheses and everything inside is called a *block*. Matched pairs can be nested. A block that is not contained within another one is called an *outer block*. For example, in each *n*-range of , the minimum element consists of a single (outer) block, and the maximum element has $\lfloor n/2 \rfloor$ outer blocks. Another example, element $w_{421} \equiv (0(0))0()$ contains two outer blocks: left block $(0(0))$ and right one $()$; and the first block also contains an inner block $(0)$. In this section, we work only with outer blocks, hereinafter just blocks.

**Definition 2.2.** Let $x$ be a nonzero word in $\mathfrak{M}$, and let's select an arbitrary (outer) block into $x$ and replace all external parentheses with zeros. As a result, we get an *extended block* of same length with leading and trailing zeros in the general case.

We assume that leading zeros only formally extend the code without changing Motzkin words themselves (as in the case of natural numbers). In other words, in the extended block, we can add or remove arbitrarily leading zeros if necessary.

In the example below, we introduce addition $\oplus$ and subtraction $\ominus$ for the Motzkin words.

**Example 2.3.** Let's choose $w_{736} \equiv ( )0(0( ))0$. In this case, we can get two extended blocks: the left block $( )0000000 \equiv w_{708}$ and the right one $000(0( ))0 = (0( ))0 \equiv w_{28}$. And immediately we get the index equality $708 + 28 = 736$. Let's repeat this equality for Motzkin words in a column:

$$\begin{array}{r} ( )0 0 0 0 0 0 0 \\ \oplus\ 0 0 0 ( 0 ( ) ) 0 \\ =\ ( ) 0 ( 0 ( ) ) 0. \end{array}$$

The leading zeros came in handy in the second block. In fact, we *decomposed* one element of $\mathfrak{M}$ into two words: $w_{736} = w_{708} \oplus w_{28}$. Also in the column, we can write down two corresponding subtraction operations: $w_{708} = w_{736} \ominus w_{28}$ and $w_{28} = w_{736} \ominus w_{708}$. □

Note that when performing arithmetic operation, zero in Motzkin word is processed like zero in integer. Let's write it in the form of the following rules:

$$0 \oplus 0 = 0, \quad 0 \oplus ( = ( \oplus 0 = (, \quad 0 \oplus ) = ) \oplus 0 = ),$$
$$0 \ominus 0 = 0, \quad ( \ominus 0 = (, \quad ) \ominus 0 = ), \quad ( \ominus ( = 0, \quad ) \ominus ) = 0.$$

Obviously, for $w_0 \equiv 0$ and any $x \in \mathfrak{M}$, it is fair:

$$x \oplus w_0 = w_0 \oplus x = x,$$
$$x \ominus w_0 = x, \quad x \ominus x = w_0.$$

Now we can state the corresponding theorem.

**Theorem 2.4.** *Let the element* $x \in \mathfrak{M}$ *be divided into k extended blocks* $x_i$, *that is,*

(2.4) $$x = x_1 \oplus x_2 \oplus \ldots \oplus x_k = \bigoplus_{1 \geq i \geq k} x_i.$$

*Then the index equation takes the form*

(2.5) $$\text{ind } x = \sum_{1 \geq i \geq k} \text{ind } x_i.$$

*Proof.* In the case $k = 1$, the theorem is obvious. Let $k > 1$, and let $x_1 > x_2 > \ldots > x_k$. Then in (2.4) the extended blocks are listed in descending order of ranges (leading zeros do not count). This means that the source element $x$ and the first block $x_1$ are the same length, i.e., are in the same range. Obviously, the first block has trailing zeros, and the size of the tail is not less than the length of the second block $x_2$.

Let's "walk" through Motzkin Row from $x_1$ to $x_2$. Step by step we will iterate over the following elements: $x_1 \oplus w_1$, $x_1 \oplus w_2$, $x_1 \oplus w_3$, and so on. With each step the index of the sum is incremented by 1. As a result, we get the element $x_1 \oplus x_2$ with index ind $x_1$ + ind $x_2$. Next, repeat the procedure with the passage of Motzkin Row from point $x_1 \oplus x_2$ to point $(x_1 \oplus x_2) \oplus x_3 = x_1 \oplus x_2 \oplus x_3$. The index new point is ind $x_1$ + ind $x_2$ + ind $x_3$. Further by induction, we obtain the final equality (2.5). □

In Example 2.3, we decomposed a word of $\mathfrak{M}$ into extended blocks. In the following example, we insert a block into the zero zone of a given element.

**Example 2.5.** The 9-word $w_{710} \equiv (\,)0000(0)$ has an inner fragment of four zeros. Let's fill this place with the 7-word $w_{72} \equiv (0)0000$. First add the indices $710 + 72 = 782$. The corresponding element is $w_{782} \equiv (\,)(0)0(0)$. Summing up the brackets, we get an identical result:

$$
\begin{aligned}
&(\,)0\,0\,0\,0\,(\,0\,) \\
\oplus\ &0\,0\,(\,0\,)\,0\,0\,0\,0 \\
=\ &(\,)\,(\,0\,)\,0\,(\,0\,).
\end{aligned}
$$

Again, for convenience, we put leading zeros in the second summand. Also we can write two subtraction operations: $w_{710} = w_{782} \ominus w_{72}$ and $w_{72} = w_{782} \ominus w_{710}$. Given element $w_{710}$ contains two blocks $(\,)0000000 \equiv w_{708}$ and $000000(0) = (0) \equiv w_2$. So we can write down the sum with three extended blocks $w_2 \oplus w_{72} \oplus w_{710} = w_{782}$. □

In the general case, any nonzero element of Motzkin Row can be decomposed into extended blocks, which can be further grouped arbitrarily.

**Definition 2.6.** The words $x, y \in \mathfrak{M}$, $x < y$, are called *noncrossing* if each outer block in $x$ corresponds in $y$ to the zero zone between the outer blocks in the same positions. In this case, addition $\oplus$ is defined for these words, i.e., $x \oplus y \in \mathfrak{M}$.

Noncrossing words have different sizes. The operation $\oplus$ has many properties of ordinary addition of numbers. This operation is *symmetric*, that is, for noncrossing $x, y \in \mathfrak{M}$ $x \oplus y = y \oplus x \in \mathfrak{M}$. Or here's another: the sum of two nonzero words is greater than each of the terms.

Here is the definition associated with the subtraction operation.

**Definition 2.7.** Let $x, y \in \mathfrak{M}$, and let $x > y$. We say that $y$ is *included* in $x$, $y \sqsubset x$, if each extended block from $y$ is also contained in $x$. In this case, subtraction $\ominus$ is defined for these words, i.e., $x \ominus y \in \mathfrak{M}$.

It's easy to see, for $x, y, z \in \mathfrak{M}$, if $x \oplus y = z$ then $z \ominus x = y$ and $z \ominus y = x$.

**Example 2.8.** Let's continue Example 2.5. We divided $w_{782} = ( )(0)0(0)$ into three extended blocks $w_2 = (0)$, $w_{72} = (0)0000$ and $w_{708} = ( )0000000$. These blocks can be grouped in various ways

$$(w_2 \oplus w_{72}) \oplus w_{708} = w_2 \oplus (w_{72} \oplus w_{708}) = (w_2 \oplus w_{708}) \oplus w_{72} = w_{782}.$$

Thus the result does not change when the extended blocks are regrouped. Let's write another equality:

$$w_2 \oplus w_{708} = w_{782} \ominus w_{72}.$$

The last equality contains both arithmetic operations on elements of $\mathfrak{M}$. □

Theorem 2.4 implies the following corollary for disjoint and nested words.

**Corollary 2.9.** *Let $x, y \in \mathfrak{M}$, and let $x > y$.*

(i) *If $x$ and $y$ do not cross, then $x \oplus y = z \in \mathfrak{M}$ and $\text{ind } z = \text{ind } x + \text{ind } y$.*

(ii) *If $y$ included in $x$, or $y \sqsubset x$, then $x \ominus y = z \in \mathfrak{M}$ and $\text{ind } z = \text{ind } x - \text{ind } y$.*

This concludes the arithmetic on the elements of $\mathfrak{M}$. Let's move on to logical operations that allow you to modify words, i.e., to perform throws on Motzkin Row.

## 3 Navigation along Motzkin Row

In the previous section, we operated with words of Motzkin Row at the level of outer and extended blocks using arithmetic operations. Now we will modify blocks to move through Motzkin Row. All transformations in each block are reflected in the same way in the Motzkin word where this block is included. For example, let $x, y \in \mathfrak{M}$, and let $y \sqsubset x$; then in accordance with Corollary 2.8, corrections in $y$ will equally affect both index $y$ and index $x$, that is, the increments of both indices are the same: $\Delta \text{ ind } y = \Delta \text{ ind } x$.

In this section, we will discuss "floating" brackets, i.e. parentheses "drift" towards adjacent zeros. Let's start with the simplest, with the left parenthesis of the outer block of a word from $\mathfrak{M}$.

3.1. **Wandering left parenthesis**. Recall that in the $n$-range of Motzkin Row, the initial $n$-word $(00\ldots 0)$ has an index of $M_{n-1}$; in the next $(n+1)$-range, the minimum word has an index of $M_n$. To go from one element to another, you need to make a jump up the row by the value of $M_n - M_{n-1}$. In this case, the left parenthesis of the Motzkin word moves one position to the left: from the $n$-th symbol to the $(n+1)$-th symbol. (The symbols are numbered from right to left as in natural numbers!) If we put the leading zero before the word (useless zero again came in handy), then we can say that the left parenthesis and the adjacent zero are reversed.

Let's explain again, in the outer block, we move the leftmost parenthesis to the left. In this case, the magnitude of the jump along Motzkin Row depends only on the position of a parenthesis in the word and the presence of an adjacent left-side zero (possibly leading zero). We can say that in the corresponding extended block the range is increased by one. The contents of the outer block may be arbitrary, and the following obvious statement about this.

**Corollary 3.1.** *Let $x \in \mathfrak{M}$ has the outer block with the leftmost parenthesis in the k-th position, and let the adjacent (or leading) zero be to the left of the parenthesis. Then a permutation of the bracket and zero increases the index of x by $M_k - M_{k-1}$, that is,*

$$\Delta \operatorname{ind} x = M_k - M_{k-1}.$$

Let's display the procedure for moving the left parenthesis as follows:

$$w_4 \equiv (00) = 0(00) \implies 0(^{+1}00) = (000) \equiv w_9.$$

We marked in red the movable parenthesis and plus one in the superscript (moving one position to the left with an increase of the word index). In this case, the jump is $M_4 - M_3 = 9 - 4 = 5$. Let's take two consecutive steps:

$$w_4 \equiv (00) \implies (^{+1}00) \implies (^{+1}000) = (0000) \equiv w_{21} \quad \text{or} \quad (^{+2}00) = (0000).$$

In this case, the jump is $(M_4 - M_3) + (M_5 - M_4) = M_5 - M_3 = 21 - 4 = 17$.

Obviously, you can move a parenthesis back from left to right (superscript negative), and then the word index decreases. The plus sign can be omitted in the superscript. A similar notation can be used for words of indefinite length. For example, the minimum element of the *n*-range of Motzkin Row can be written as follows:

$$\min \mathfrak{M}_n = (^{n-2}), \quad n \geq 2.$$

It is easy to see when the left bracket from position $k$ is shifted by $\pm j$ positions (left $+j$ or right $-j$), the index increment is

(3.1) $$\delta_k^{\,j} = M_{k-1\pm j} - M_{k-1}, \quad j \geq 0.$$

Consider the example.

**Example 3.2.** Let's fix the movable left bracket in the 6th position ($k = 6$). Then, according to (3.1), the increment of the index is $M_{5+j} - M_5$. Below, for the selected Motzkin words, we obtain the index increment for cases $k = \pm 1$ and $k = \pm 2$.

$w_{21} \equiv (0000) \implies (^{+1}0000) = (00000) \equiv w_{51};\ \Delta \operatorname{ind} w_{21} = M_6 - M_5 = 51 - 21 = +30.$
$w_{28} \equiv (0())0 \implies (^{-1}0())0 = (())0 \equiv w_{16};\quad \Delta \operatorname{ind} w_{28} = M_4 - M_5 = 9 - 21 = -12.$
$w_{50} \equiv ()()() \implies (^{+2})()() = (00)()() \equiv w_{156};\ \Delta \operatorname{ind} w_{50} = M_7 - M_5 = 127 - 21 = +106.$
$w_{294} \equiv ()(000)0 \implies ()(^{-2}000)0 = ()00(0)0 \equiv w_{277};\ \Delta \operatorname{ind} w_{294} = M_3 - M_5 = 4 - 21 = -17.$
$w_{742} \equiv ()0((0)0) \implies ()0(^{+1}(0)0) = ()(0(0)0) \equiv w_{772};\ \Delta \operatorname{ind} w_{742} = M_6 - M_5 = 51 - 21 = +30.$

Pay attention to the various elements of Motzkin Row, the index increments coincide if $k$ and $j$ are fixed (see $w_{21}$ and $w_{742}$). □

3.2. **Control points in Motzkin row.** In Motzkin Row, we can consider the minimum/maximum elements as *control points* (landmarks, pegs) since their indices are known. Control points allow you to speed up the identification procedure: calculate the index of a given Motzkin word, or vice versa, construct an element of Motzkin Row using a given index. Our analysis allows us to obtain additional landmarks.

(1) Let $a^{(n)}$ be the *n*-word from $\mathfrak{M}$, and let's define $a^{(n)}$ as follow:

$$a^{(n)} = \max \mathfrak{M}_n \ominus \max \mathfrak{M}_{n-2} = (_n)()\ldots \ominus (_{n-2})()\ldots = (_n) 0\ldots 0.$$

The subscript next to the bracket is the bracket position in the word. Obviously,

$$\text{ind } a^{(n)} = (M_n - 1) - (M_{n-2} - 1) = M_n - M_{n-2}.$$

(2) Let's calculate another landmark in the *n*-range: $b^{(n)} = (_n 0) 0\ldots 0$. Take the intermediate (*n*–1)-word

$$a^{(n-1)} = (_{n-1}) 0\ldots 0 = \max \mathfrak{M}_{n-1} \ominus \max \mathfrak{M}_{n-3}, \quad \text{ind } a^{(n-1)} = M_{n-1} - M_{n-3}.$$

In $a^{(n-1)}$, we move the left bracket by one position to the left. As a result, we get the following index polynomial (see Corollary 3.1):

$$\text{ind } b^{(n)} = (M_{n-1} - M_{n-3}) + (M_{n-1} - M_{n-2}) = 2M_{n-1} - M_{n-2} - M_{n-3}.$$

(3) Here is one more landmark in the *n*-range:

$$c^{(n)} = (_n 0) ()()\ldots = b^{(n)} \oplus \max \mathfrak{M}_{n-3}, \quad \text{ind } c^{(n)} = \text{ind } b^{(n)} + M_{n-3} - 1.$$

Thus, for the *n*-range of Motzkin Row, there are different landmarks for solving the identification problem. In Table 1 below, we have listed the seven control points of the *n*-range in order of increasing indices. Five points we have considered, two points $d^{(n)}$ and $h^{(n)}$ will be considered further (so their indices are not specified). If desired, the reader can independently calculate additional landmarks.

Table 1.

| 1 | min $\mathfrak{M}_n$ | $= (00\ldots 0)$ | ind min $\mathfrak{M}_n$ | $= M_{n-1}$ |
|---|---|---|---|---|
| 2 | $d^{(n)}$ | $= (0( ))00\ldots 0$ | | |
| 3 | $b^{(n)}$ | $= (0)00\ldots 0$ | ind $b^{(n)}$ | $= 2M_{n-1} - M_{n-2} - M_{n-3}$ |
| 4 | $c^{(n)}$ | $= (0)()( )\ldots ( )[0]$ | ind $c^{(n)}$ | $= 2M_{n-1} - M_{n-2} - 1$ |
| 5 | $h^{(n)}$ | $= ((0))00\ldots 0$ | | |
| 6 | $a^{(n)}$ | $= ( )00\ldots 0$ | ind $a^{(n)}$ | $= M_n - M_{n-2}$ |
| 7 | max $\mathfrak{M}_n$ | $= ( )( )\ldots ( )[0]$ | ind max $\mathfrak{M}_n$ | $= M_n - 1$ |

3.3. **Wandering right parenthesis.** As in the case with left brackets, we work only with external right brackets of blocks (both outer and extended). With the right brackets, everything is not so obvious, and yet it is not difficult to obtain a formula similar to (3.1). As before, we will record the change in the word index when the parenthesis and adjacent zero are reversed. Note that, unlike left brackets, moving the right parentheses never changes the element range.

Let's choose a simple word from $\mathfrak{M}$ with the right parenthesis in the $k$-th position and with zero in the $(k+1)$-th position. One of the control points discussed above will suit us:

$$b^{(k+2)} = (0)_k 00\ldots 0, \quad \text{ind } b^{(k+2)} = 2M_{k+1} - M_k - M_{k-1}.$$

At the tail of this $(k+2)$-word, there are $k-1$ zeros. We again marked the right parenthesis by its number. The position of the left parenthesis is $k+2$. After the permutation of the right parenthesis with left zero, we get another control point (see Table 1):

$$a^{(k+2)} = (\ )_{k+1} 000\ldots 0, \quad \text{ind } a^{(k+2)} = M_{k+2} - M_k.$$

The position of the right bracket is increased by one, and in the tail of $a^{(k+2)}$ there are $k$ zeros. The index of $b^{(k+2)}$ has changed as follows:

$$\begin{aligned}\Delta \text{ind } b^{(k+2)} &= \text{ind } a^{(k+2)} - \text{ind } b^{(k+2)} \\ &= (M_{k+2} - M_k) - (2M_{k+1} - M_k - M_{k-1}) \\ &= M_{k+2} - 2M_{k+1} + M_{k-1}, \quad k \geq 1.\end{aligned}$$

In this case, the selected control point $b^{(k+2)}$ is the simple extended block. It is easy to see that the resulting index polynomial is valid for more complex blocks.

**Corollary 3.3.** *Let $x \in \mathfrak{M}$ has the outer block with the rightmost parenthesis in the $k$-th position, and let the adjacent zero be to the left of the parenthesis. Then a permutation of the bracket and zero increases the index of $x$ by*

(3.2) $$\xi_k = M_{k+2} - 2M_{k+1} + M_{k-1}, \quad k \geq 1.$$

Naturally, the right bracket can drift to the right (if there is adjacent zero on the right). Then in the index polynomial (3.2), it is enough to change the signs and reduce $k$ by one (to the right of the equal sign). Let us consider the following example.

**Example 3.4.** Let's choose the words with a movable right parenthesis in the same position, and let $k = 5$. We assume that the selected right parenthesis drifts to the left. Then $\xi_5 = M_7 - 2M_6 + M_4 = 127 - 2\times 51 + 9 = 34$.

$w_{72} \equiv (0)0000 \Rightarrow (\ )00000 \equiv w_{106}$; $\Delta \text{ind } w_{72} = 106 - 72 = 34$;
$w_{154} \equiv (00)(()) \Rightarrow (0)0(()) \equiv w_{188}$; $\Delta \text{ind } w_{154} = 188 - 154 = 34$;
$w_{658} \equiv (()0)(0)0 \Rightarrow (())0(0)0 \equiv w_{692}$; $\Delta \text{ind } w_{658} = 692 - 658 = 34$.

We highlighted in red the floating parenthesis and the left adjacent zero. □

The numbers $\xi_k$ obtained by formula (3.2) for $k = 1, 2, 3, ...$ form the sequence:

(3.3)    1, 2, 5, 13, 34, 90, 240, 645, 1745, 4750, 13001, 35762, 98815, 274158, …

Let's summarize. In any word of Motzkin Row, we fix the *n*-th position of the rightmost bracket of an outer block; also zero is required to the left of the bracket. If the right parenthesis is swapped with zero, then the word index increases in accordance with the index polynomial (3.2).

## 4    Merge and split blocks

We examined simple operations on the elements of Motzkin Row. Maybe over time this lexicographic sequence will become popular, and then the study of the properties of Motzkin Row will be done by various mathematicians.

Previously, we grouped outer blocks and moved their borders. In these cases, the range of words sometimes changed. Below we consider a few additional properties that allow you to change the number of blocks in the Motzkin word. By modifying the outer block inside, you can split it into several blocks. Accordingly, reverse procedures work. In the properties below, the range of words does not change. So, let's look inside the block.

4.1. **Remove parentheses / insert parentheses**. Consider at the procedure for removing the extreme brackets of neighboring blocks, and the reverse procedure for restoring the brackets. First, consider a simple word with two blocks that are separated by zeros (see the index polynomials in Table 1):

$$x = (\ )_l\, 0...0\, (_k\, 0...0) = a^{(l+1)} \oplus \min \mathfrak{M}_k, \quad \text{ind}\, x = (M_{l+1} - M_{l-1}) + M_{k-1}.$$

In the ($l$+1)-word $x$, we have marked red brackets for deletion, namely: the last parenthesis of the left block in the *l*-th position and the initial parenthesis of the right block in the *k*-th position, $l > k$ (recall, the brackets are numbered from right to left, from the end of the word, as in integers).

Replace the red parentheses with zeros. As a result, we get the word $\min \mathfrak{M}_{l+1}$ with the index $M_l$. It is obvious that the index of $x$ is reduced by

$$-\Delta\, \text{ind}\, x\, =\, (M_{l+1} - M_{l-1} + M_{k-1})\, - M_l.$$

Let's complicate the given word a bit (again, see Table 1):

$$x' = (0)_l\, 0...0\, (_k\, 0...0) = b^{(l+2)} \oplus \min \mathfrak{M}_k, \quad \text{ind}\, x' = (2M_{l+1} - M_l - M_{l-1}) + M_{k-1}.$$

After removing the red brackets, we get $\min \mathfrak{M}_{l+2}$ with the index $M_{l+1}$. In this case, the index of $x'$ is reduced by

$$-\Delta\, \text{ind}\, x'\, =\, (2M_{l+1} - M_l - M_{l-1} + M_{k-1})\, - M_{l+1}.$$

As you can see, the same index polynomial is obtained. You can arbitrarily complicate the contents of neighboring blocks, add additional blocks in the word $x$ both at the beginning and at the end, but the index polynomial will be the same. The result depends only on the location of the brackets to be removed. Let us formulate an appropriate statement.

**Corollary 4.1.** *In a word of Motzkin Row, let the neighboring outer blocks be separated by a zero zone in the boundaries of $l$ and $k$, $l > k$. That is, the outer parentheses of the blocks are in the specified points. Then replacing these two brackets with zeros reduces the word index by*

(4.1) $\qquad \zeta_{k,l} = M_{l+1} - M_l - M_{l-1} + M_{k-1}, \ l > k \geq 2.$

Obviously, the index polynomial (4.1) works both ways. We can not only delete the extreme brackets of neighboring blocks, but also perform the reverse procedure. You can insert a pair of oppositely directed parentheses in the zero zone of the selected block.

**Example 4.2.** In the words below, modifiable characters are shown in red. Recall the first ten Motzkin numbers $M_n$, $n = 1, 2, \ldots$ : 1, 2, 4, 9, 21, 51, 127, 323, 835, 2188.
$w_{278} \equiv ()00(()) \Leftrightarrow (0000()) \equiv w_{129}$, $\zeta_{4,7} = M_8 - M_7 - M_6 + M_3 = 323 - 127 - 51 + 4 = 149$.
$w_{491} \equiv (0)(00()) \Leftrightarrow (0)()(()) \equiv w_{516}$, $\zeta_{4,5} = M_6 - M_5 - M_4 + M_3 = 51 - 21 - 9 + 4 = 25$.
$w_{1152} \equiv (0()00)000 \Leftrightarrow (0())()000 \equiv w_{1216}$, $\zeta_{5,6} = M_7 - M_6 - M_5 + M_4 = 127 - 51 - 21 + 9 = 64$.
$w_{2153} \equiv ()()0(())0 \Leftrightarrow ()(000())0 \equiv w_{1999}$, $\zeta_{5,7} = M_8 - M_7 - M_6 + M_4 = 323 - 127 - 51 + 9 = 154$.
□

Adjacent blocks are often processed, this is when $l = k + 1$. In this case

(4.2) $\qquad \zeta_{k,k+1} = M_{k+2} - M_{k+1} - M_k + M_{k-1}, \ k \geq 2.$

The numbers $\zeta_{k,k+1}$ obtained by (4.2) for $k = 2, 3, \ldots$ form the sequence:

(4.3) $\quad$ 4, 10, 25, 64, 166, 436, 1157, 3098, 8360, 22714, 62086, 170614, …

4.2. **Permutation of adjacent parentheses**. Let $x \in \mathfrak{M}$, and let $x$ contains two adjacent outer blocks:

$$x = \ldots ( \ldots )^{\leftrightarrow} (_k \ldots ) \ldots$$

We have shown the outer parentheses of the blocks. Dots mark arbitrary fragments of $x$. The touching brackets (positions $k$ and $k+1$) are highlighted in red.

Let's permute the marked brackets. In this case, two blocks are *merged*:

$$y = \ldots ( \ldots ()_k \ldots ) \ldots$$

Then the index $y$ is calculated simply:

(4.4) $\qquad \text{ind } y = \text{ind } x - M_k \ \text{ or } \ \Delta \text{ ind } x = -M_k.$

There is no strict proof of the property (4.4) yet. Let's formulate a guess.

**Conjecture 4.3** (merging of blocks by reorienting brackets). *Let $x \in \mathfrak{M}$, and let $x$ has two adjacent outer blocks that contact at the positions $k$ and $k+1$. Then a permutation of the contact brackets reduces the index of $x$ by $M_k$, i.e., $\Delta \operatorname{ind} x = -M_k$.*

Obviously, the reverse transition from $y$ to $x$ leads to *splitting* of the block. The author tested Conjecture 4.1 in the first ten ranges of Motzkin Row.

In this regard, we consider example, in which an additional checkpoint between $\min \mathfrak{M}_n$ and $b^{(n)}$ is calculated (see Table 1 in section 3.2).

**Example 4.4.** Let's sum up two noncrossing words $b^{(n)}$ and $a^{(n-3)}$:

$$r^{(n)} = b^{(n)} \oplus a^{(n-3)} = (_n 0) 00\ldots 0 \oplus (_{n-3}) 00\ldots 0 = (_n 0) (_{n-3}) 00\ldots 0.$$

In the case of permutation of two red brackets, we obtain a new control point

$$d^{(n)} = (0(\;))00\ldots 0, \quad \operatorname{ind} d^{(n)} = \operatorname{ind} r^{(n)} - M_{n-3}.$$

Let's write the result below into the Table 1.

$$\begin{aligned}\operatorname{ind} d^{(n)} &= (\operatorname{ind} b^{(n)} + \operatorname{ind} a^{(n-3)}) - M_{n-3} \\ &= (2M_{n-1} - M_{n-2} - M_{n-3}) + (M_{n-3} - M_{n-5}) - M_{n-3} \\ &= 2M_{n-1} - M_{n-2} - M_{n-3} - M_{n-5}. \qquad \square\end{aligned}$$

**4.3. Permutation of non-adjacent parentheses**. Let $x \in \mathfrak{M}$, and let $x$ contains two outer blocks with zero between them:

$$x = \ldots (\ldots) 0 (_k \ldots) \ldots$$

Again dots mark arbitrary fragments of $x$. The moving brackets are shown in red; in the $(k+1)$-position, there is a zero between the brackets.

As a result of the permutation, we get the word with a smaller index:

$$y = \ldots (\ldots (0)_k \ldots) \ldots$$

Let's denote such a change in index $\psi_k = \operatorname{ind} x - \operatorname{ind} y$.

The formula of $\psi_k$ is unknown to the author. The analysis of the available elements of $\mathfrak{M}$ gives such a sequence of numbers for $\psi_k$, $k = 2, 3, \ldots$:

(4.5) $\qquad$ 4, 10, 25, 65, 171, 456, 1227, 3328, 9084, …

In the example below, using (4.5), we move the bracket inside the block to get the control point in the $n$-range between $c^{(n)}$ and $a^{(n)}$ (see Table 1 in section 3.2).

**Example 4.5.** In the $n$-word $d^{(n)} = (0(\;))00\ldots 0$ (see Example 4.4), let's rearrange the zero and the left parenthesis (marked in red) to get the control point $h^{(n)} = ((0))00\ldots 0$. Let's follow some steps.

First (i) divide the block into two parts: $(0(\;)) \Rightarrow (0)(\;)$ (see Conjecture 4.1, $k = n-3$). As a result we obtain the intermediate point (see Example 4.4):

$$r^{(n)} = (0)(\,)00\ldots 0, \quad \text{ind } r^{(n)} = \text{ind } d^{(n)} + M_{n-3} = 2M_{n-1} - M_{n-2} - M_{n-5}.$$

Then (ii) we move the right parenthesis to the left (see Corollary 3.3, $k = n-2$):

$$(0)(\,)00\ldots 0 \Rightarrow (\,)0(\,)00\ldots 0 = s^{(n)},$$

$$\text{ind } s^{(n)} = \text{ind } r^{(n)} + M_n - 2M_{n-1} + M_{n-3} = M_n - M_{n-2} + M_{n-3} - M_{n-5}.$$

In the last step (iii), combine the two blocks with zero between them:

$$(\,)0(_{n-3})00\ldots 0 \Rightarrow ((0))00\ldots 0 = h^{(n)}, \quad \text{ind } h^{(n)} = \text{ind } s^{(n)} - \psi_{n-3}.$$

In Table 1 we can write the following index polynomial:

$$\text{ind } h^{(n)} = M_n - M_{n-2} + M_{n-3} - M_{n-5} - \psi_{n-3}.$$

The index increment is

$$\Delta \text{ ind } d^{(n)} = \text{ind } h^{(n)} - d^{(n)} = M_n - 2M_{n-1} + 2M_{n-3} - \psi_{n-3}.$$

Check the result:

(1) $w_{70} \equiv (0(\,))00 \Rightarrow ((0))00 \equiv w_{88}$, $n = 7$,

$\Delta \text{ ind } w_{70} = M_7 - 2M_6 + 2M_4 - \psi_4 = 127 - 2\cdot 51 + 2\cdot 9 - 25 = 18$.

(2) $w_{464} \equiv (0(\,))(0)0 \Rightarrow ((0))(0)0 \equiv w_{584}$, $n = 9$,

$\Delta \text{ ind } w_{464} = M_9 - 2M_8 + 2M_6 - \psi_6 = 835 - 2\cdot 323 + 2\cdot 51 - 171 = 120$.

(3) $w_{1502} \equiv ((0)_k((\,)))0 \Rightarrow (\,)0(((\,)))0 \equiv w_{1958}$, $k = 7$. Here, only the block is divided into two pieces, between which there is a single zero. The index of $w_{1502}$ is just increased by $\psi_7 = 456$. □

Gzhel State University, Moscow, 140155, Russia
http://www.en.art-gzhel.ru/

**Addendum**. Motzkin Row (ranges 1 – 9).

000: 0, ( ), (0), ( )0, (00), (0)0, (( )), ( )00, ( )( ), (000), (00)0, (0( )), (0)00, (0)( ), ((0))
015: (( )0), (( ))0, ( )000, ( )0( ), ( )(0), ( )( )0, (0000), (000)0, (00( )), (00)00, (00)( ), (0(0))
027: (0( )0), (0( ))0, (0)000, (0)0( ), (0)(0), (0)( )0, ((00)), ((0)0), ((0))0, ((( ))), (( )00), (( )0)0
039: (( )( )), (( ))00, (( ))( ), ( )0000, ( )00( ), ( )0(0), ( )0( )0, ( )(00), ( )(0)0, ( )(( )), ( )( )00, ( )( )( )
051: (00000), (0000)0, (000( )), (000)00, (000)( ), (00(0)), (00)0)0, (00))0, (00)000
060: (00)0( ), (00)(0), (00)( )0, (0(00)), (0(0)0), (0(0))0, (0(( ))), (0( )00), (0( )0)0, (0( )( ))
070: (0( ))00, (0( ))( ), (0)0000, (0)00( ), (0)0(0), (0)0( )0, (0)(00), (0)(0)0, (0)(( )), (0)( )00
080: (0)( )( ), ((000)), ((00)0), ((00))0, ((0( ))), ((0)00), ((0)0)0, ((0)( )), ((0))00, ((0))( )
090: (((0))), (((0))), ((( )0), ((( ))0, (( )000, (( )00)0, (( )0( )), (( )0)00, (( )0)( ), (( )(0)), (( )( )0)
101: (( )( ))0, (( ))000, (( ))0( ), (( ))(0), (( ))( )0, ( )00000, ( )0000( ), ( )000(0), ( )000( )0, ( )00(00)
111: ( )0(0)0, ( )0(( )), ( )0( )00, ( )0( )( ), ( )(000), ( )(00)0, ( )(0( )), ( )(0)00, ( )(0)( ), ( )((0))
121: ( )(( )0), ( )(( ))0, ( )( )000, ( )( )0( ), ( )( )(0), ( )( )( )0, (000000), (00000)0, (0000( )), (0000)00
131: (0000)( ), (000(0)), (000(0)), (000( )0), (000))0, (000)000, (000)0( ), (000)(0), (000)( )0, (00(00))
140: (00(0)0), (00(0))0, (00(( ))), (00( )00), (00( )0)0, (00( )( )), (00( ))00, (00( ))( ), (00)0000
149: (00)00( ), (00)0(0), (00)0( )0, (00)(00), (00)(0)0, (00)(( )), (00)( )00, (00)( )( ), (0(000))
158: (0(00)0), (0(00))0, (0(0( ))), (0(0)00), (0(0)0)0, (0(0)( )), (0(0))00, (0(0))( ), (0(( )))
167: (0(( )0)), (0(( ))0), (0(( ))0)0, (0( )000), (0( )00)0, (0( )0( )), (0( )0)00, (0( )0)( ), (0( )(0))
176: (0( )( )0), (0( )( ))0, (0( ))000, (0( ))0( ), (0( ))(0), (0( ))( )0, (0)00000, (0)0000( ), (0)000(0)
185: (0)000( )0, (0)00(00), (0)00(0)0, (0)00(( )), (0)00( )00, (0)00( )( ), (0)0(000), (0)0(00)0, (0)0(0( ))
194: (0)(0)00, (0)(0)( ), (0)((0)), (0)(( )0), (0)(( ))0, (0)( )000, (0)( )0( ), (0)( )(0), (0)( )( )0
203: ((0000)), ((000)0), ((000))0, ((00( ))), ((00)00), ((00)0)0, ((00)( )), ((00))00, ((00))( )
212: ((0(0))), ((0( )0)), ((0( ))0), ((0( )))0, ((0)000, ((0)00)0, ((0)0( )), ((0)0)00, ((0)0)( )
221: ((0)(0)), ((0)( )0), ((0)( ))0, ((0))000, ((0))0( ), ((0))(0), ((0))( )0, (((00))), (((0)0))
230: (((0))0), ((( ))0)0, ((( ))( )), ((( )00), ((( )0)0), ((( )0)0, ((( )( )), ((( ))00, ((( ))0)0
239: (((( )))), ((( ))00, ((( ))( ), (( )0000), (( )000)0, (( )00( )), (( )00)00, (( )00)( ), (( )0(0))
248: (( )0( )0), (( )0( ))0, (( )0)000, (( )0)0( ), (( )0)(0), (( )0)( )0, (( )(00)), (( )(0)0), (( )(0))0
257: (( )(( ))), (( )( )00), (( )( )0)0, (( )( )( )), (( )( ))00, (( )( ))( ), (( ))0000, (( ))00( ), (( ))0(0)
266: (( ))0( )0, (( ))(00), (( ))(0)0, (( ))(( )), (( ))( )00, (( ))( )( ), ( )000000, ( )00000( ), ( )0000(0)
275: ( )000( )0, ( )000(00), ( )000(0)0, ( )00(( )), ( )000( )0, ( )000( )( ), ( )00(000), ( )00(00)0, ( )00(0( ))
284: ( )0(0)00, ( )0(0)( ), ( )0((0)), ( )0(( )0), ( )0(( ))0, ( )00( )000, ( )0(0)( ), ( )00( )(0), ( )0(0)( )0
293: ( )(0000), ( )(000)0, ( )(00( )), ( )(00)00, ( )(00)( ), ( )(0(0)), ( )(0( )0), ( )(0( ))0, ( )(0)000
302: ( )(0)0( ), ( )(0)(0), ( )(0)( )0, ( )((00)), ( )((0)0), ( )((0))0, ( )((( ))), ( )(( )00), ( )(( )0)0
311: ( )(( )( )), ( )(( ))00, ( )(( ))( ), ( )( )0000, ( )( )00( ), ( )( )0(0), ( )( )0( )0, ( )( )(00), ( )( )(0)0
320: ( )( )(( )), ( )( )( )00, ( )( )( )( ), (0000000), (000000)0, (00000( )), (00000)00, (00000)( )
328: (0000(0)), (0000(0)), (0000( ))0, (0000)000, (0000)0( ), (0000)(0), (0000)( )0, (000(00))
336: (000(0)0), (000(0))0, (000(( ))), (000( )00), (000( )0)0, (000( )( )), (000( ))00, (000( ))( )
344: (000)0000, (000)00( ), (000)0(0), (000)0( )0, (000)(00), (000)(0)0, (000)(( )), (000)( )00
352: (000)( )( ), (00(000)), (00(00)0), (00(00))0, (00(0( ))), (00(0)00), (00(0)0)0, (00(0)( ))
360: (00(0))00, (00(0))( ), (00(( ))), (00(( )0), (00(( ))0), (00(( ))0)0, (00( )000), (00( )00)0
368: (00( )0( )), (00( )0)0, (00( )0)( ), (00( )(0)), (00( )( )0), (00( )( ))0, (00( ))000, (00( ))0( )
376: (00( ))(0), (00( ))( )0, (00)00000, (00)0000( ), (00)000(0), (00)000( )0, (00)00(00), (00)00(0)0
384: (00)0(( )), (00)0( )00, (00)0( )( ), (00)(000), (00)(00)0, (00)(0( )), (00)(0)00, (00)(0)( )
392: (00)((0)), (00)(( )0), (00)(( ))0, (00)( )000, (00)( )0( ), (00)( )(0), (00)( )( )0, (0(0000))
400: (0(000)0), (0(000))0, (0(00( ))), (0(00)00), (0(00)0)0, (0(00)( )), (0(00))00, (0(00))( )
408: (0(0(0))), (0(0( )0)), (0(0( ))0), (0(0( )))0, (0(0)000), (0(0)00)0, (0(0)0( )), (0(0)0)00

416: (0(0)0)(), (0(0)(0)), (0(0))00, (0(0)())0, (0(0))000, (0(0))0(), (0(0))(0), (0(0))()0
424: (0((00))), (0((0)0)), (0((0))0), (0((0)))0, (0((())), (0(()00)), (0(()0)0), (0(()0))0
432: (0(()())), (0(())00), (0(())0)0, (0(())()), (0(()))00, (0(()))(), (0()0000), (0()000)0
440: (0()00()), (0()0)00, (0()00)(), (0()0(0)), (0()0()0), (0()0))0, (0()0)000, (0()0)0()
448: (0()0)(0), (0()0)()0, (0()(00)), (0()(0)0), (0()(0))0, (0()(())), (0()()00), (0()()0)0
456: (0()()()), (0()())00, (0()())()·, (0())0000, (0())00(), (0())0(0), (0())0()0, (0())(00)
464: (0())(0)0, (0())(()), (0()))00, (0()))()·, (0)000000, (0)0000()·, (0)000(0), (0)000()0
472: (0)00(00), (0)00(0)0, (0)00(()), (0)00()00, (0)00()()·, (0)0(000), (0)0(00)0, (0)0(0())
480: (0)0(0)00, (0)0(0)()·, (0)0((0)), (0)0(()0), (0)0(())0, (0)0()000, (0)0()0()·, (0)0()(0)
488: (0)0()()0, (0)(0000), (0)(000)0, (0)(00()), (0)(00)00, (0)(00)()·, (0)(0(0)), (0)(0()0)
496: (0)(0())0, (0)(0)000, (0)(0)0()·, (0)(0)(0), (0)(0)()0, (0)((00)), (0)((0)0), (0)((0))0
504: (0)((())), (0)(()00), (0)(()0)0, (0)(()()), (0)(())00, (0)(())()·, (0)()0000, (0)()000()·
512: (0)()0(0), (0)()0()0, (0)()(00), (0)()(0)0, (0)()(()), (0)()()00, (0)()()0)·, ((00000))
520: ((0000)0), ((0000))0, ((000())), ((000)00), ((000)0)0, ((000)()), ((000))00, ((000))()·
528: ((00(0))), ((00()0)), ((00())0), ((00()))0, ((00)000), ((00)00)0, ((00)0()), ((00)0)00
536: ((00)0)(), ((00)(0)), ((00)(0)·), ((00)())0, ((00))000, ((00))0()·, ((00))(0), ((00))()0
544: ((0(00))), ((0(0)0)), ((0(0))0), ((0(0)))0, ((0(())), ((0()00)), ((0()0)0), ((0()0))0
552: ((0()())), ((0())00), ((0())0)0, ((0())()), ((0()))00, ((0()))()·, ((0)0000), ((0)000)0
560: ((0)00()), ((0)00)00, ((0)00)()·, ((0)0(0)), ((0)0()0), ((0)0))0, ((0)0)000, ((0)0)0()·
568: ((0)0)(0), ((0)0)()0, ((0)(00)), ((0)(0)0), ((0)(0))0, ((0)(())), ((0)()00), ((0)()0)0
576: ((0)()()), ((0)())00, ((0)())()·, ((0))0000, ((0))00()·, ((0))0(0), ((0))0()0, ((0))(00)
584: ((0))(0)0, ((0))(()), ((0)))00, ((0)))()·, (((000))), (((00)0)), (((00))0), (((00)))0
592: (((0())), (((0)00)), (((0)0)0), (((0)0))0, (((0)())), (((0))00), (((0))0)0, (((0))())
600: (((0)))00, (((0)))()·, ((((0))), ((((0))), ((((()))), (((()))0, (((())·)0, ((()000))
608: ((()00)0), ((()00))0, ((()0())), ((()0)00), ((()0)0)0, ((()0)()), ((()0))00, ((()0))()·
616: ((()(0))), ((()()0)), ((()()0), ((()()))0, ((())000), ((())00)0, ((())0()), ((())0)00
624: ((())0)(), ((()))(0)), ((())()0), ((())())0, ((()))000, ((()))0()·, ((()))·(0), ((()))·()0
632: (()00000), (()0000)0, (()000()), (()000)00, (()000)(), (()00(0)), (()00()0), (()00())0
640: (()00)000, (()00)0()·, (()00)(0), (()00)()0, (()0(00)), (()0(0)0), (()0(0))0, (()0(())·)
648: (()0()00), (()0()0)0, (()0()()), (()0())00, (()0())()·, (()0)0000, (()0)000()·, (()0)0(0)
656: (()0)0()0, (()0)(00), (()0)(0)0, (()0)(()), (()0)()00, (()0)()0)·, (()0))000, (()0))0)·
664: (()(00))0, (()(0())), (()(0)00), (()(0)0)0, (()(0)()), (()(0))00, (()(0))()·, (()((0)))
672: (()(()0)), (()(())0), (()(()))0, (()()000), (()()00)0, (()()0()), (()()0)00, (()()0)()·
680: (()()(0)), (()()()0), (()()())0, (()())000, (()())00)·, (()())(0), (()())()0, (())00000
688: (())000()·, (())00(0), (())00()0, (())0(00), (())0(0)0, (())0(()), (())0()00, (())0()0)·
696: (())(000), (())(00)0, (())(0()), (())(0)00, (())(0)()·, (())((0)), (())(()0), (())(())0
704: (())()000, (())()00)·, (())()0(0), (())()0)·0, ()0000000, ()000000()·, ()00000(0), ()00000()0
712: ()0000(00), ()0000(0)0, ()0000(()), ()0000()00, ()0000()()·, ()000(000), ()000(00)0, ()000(0())
720: ()000(0)00, ()000(0)()·, ()000((0)), ()000(()0), ()000(())0, ()000()000, ()000()0()·, ()000()·(0)
728: ()000()()0, ()00(0000), ()00(000)0, ()00(00()), ()00(00)00, ()00(00)()·, ()00(0(0)), ()00(0()0)
736: ()00(0())0, ()00(0)000, ()00(0)0()·, ()00(0)(0), ()00(0)()0, ()00((00)), ()00((0)0), ()00((0))0
744: ()00((()·)), ()00(()00), ()00(()0)0, ()00(()()), ()00(())00, ()00(())()·, ()00()0000, ()00()000()·
752: ()00()0(0), ()00()0()·0, ()00()(00), ()00()(0)0, ()00()(()), ()00()()00, ()00()()0)·, ()0(00000)
760: ()0(0000)0, ()0(000()), ()0(000)00, ()0(000)(), ()0(00(0)), ()0(00()0), ()0(00())0, ()0(00)000
768: ()0(00)0()·, ()0(00)(0), ()0(00)()0, ()0(0(00)), ()0(0(0)0), ()0(0(0))0, ()0(0(())·), ()0(0()00)
776: ()0(0()0)0, ()0(0()()), ()0(0())00, ()0(0())()·, ()0(0)0000, ()0(0)000)·, ()0(0)0(0), ()0(0)0()0
784: ()0(0)(00), ()0(0)(0)0, ()0(0)(()·), ()0(0)()00, ()0(0)()·)·, ()0((000)), ()0((00)0), ()0((00))0

17